\documentclass[a4paper]{amsart}

\usepackage{graphicx}
\usepackage{amssymb,amsmath}
\usepackage{mathrsfs,bm}
\usepackage{tikz}
\usepackage{subfigure}
\usepackage{epstopdf}

\allowdisplaybreaks

\newcommand{\ugec}{%
\mathrel{\reflectbox{\rotatebox[origin=c]{270}{$\geq$}}}}

\theoremstyle{plain}
\newtheorem{teo}{Theorem}[section]      
\newtheorem{lem}[teo]{Lemma}
    
\newtheorem{cor}{Corollary}       
\theoremstyle{definition}               
\newtheorem{defin}{Definition}[section]
\newtheorem{ese}{Example}[section]      
\theoremstyle{remark}                   
\newtheorem{rem}{Remark}

\title[Full binary trees and Narayana numbers]{Full binary trees, Narayana numbers and two-dimensional decompositions of integers}

\author[G. Aletti]{Giacomo Aletti}
\address{Department of Mathematics, Universit\`a degli Studi di Milano, 20131 Milano, Italy}
\email{giacomo.aletti@unimi.it}

\author[A. Daziario]{Antonio Daziario}
\address{Department of Mathematics, Universit\`a degli Studi di Milano, 20131 Milano, Italy}
\email{antonio.daziario@unimi.it}



\begin{document}

\begin{abstract}
We present two models of multitype Galton-Watson trees, that we call \textit{full binary trees} and \textit{full binary trees with survivals}. We show relevant relations between these trees and the \textit{Narayana numbers} and the \textit{two-dimensional decompositions of integers}. We prove further statistical results on our models concerning the related \textit{contour processes} and \textit{offspring distributions}.  
\end{abstract}

\nocite{*}

\maketitle

\section*{Introduction}
In this paper we present novel results and properties of two special models of the multitype Galton-Watson trees, that we call \textit{full binary trees} and \textit{full binary trees with survivals}. In particular we present interesting relations between these kind of trees and the Narayana numbers and the two-dimensional decompositions of integers.

The random trees are recently used for statistical mechanics and mathematical physics models, for instance in \cite{K18}, \cite{K22} a ferromagnetic model is studied on locally tree-like random graphs. For more details on the random graphs see \cite{K21}. In \cite{K19} and \cite{K20} are studied broadcasting problems on random trees and there are found relations with some relevant mathematical physics topics. Regarding the multitype Galton-Watson trees, they have been used to achieve the Dawson-Watanabe superprocesses, that are called more simply superprocesses (see for more details \cite{K10}). 

A multitype Galton-Watson tree is characterized in law by the offspring distribution $\bm{\mu}$ and it is characterized in realization by the \textit{contour process}, $\mathscr{C}_{\tau}\left(s\right)$ (see for instance \cite{K2}). Particular importance is given to these two characteristics in this paper. Recall that the period of $\mathscr{C}_{\tau}\left(s\right)$ is equal to $2\Vert\tau\Vert$, where $\Vert\tau\Vert$ is the number of edges of $\tau$. One important tool to be considered is the moment-generating function of the period of the contour process. In a multitype Galton-Watson tree starting witha type $i$ vertex, for every $i=1,\ldots,r,$, we have proved that the moment-generating function $F_{i}\left(s\right)$ of $2\Vert\tau_{i}\Vert$ is given by the following
\begin{equation}\label{eq:relazione_durata_con_distribuzione}
F_{i}(s)=\mathbb{E}[e^{2\Vert\tau_{i}\Vert\cdot s}]=\sum_{\bm{\alpha}\in\mathbb{N}^{r}}\Big(\mu^{(i)}(\bm{\alpha})\cdot e^{2s|\bm{\alpha}|}\cdot\prod_{k=1}^{r}F_{k}(s)^{\alpha_{k}}\Big),
\end{equation}
where $r$ is the number of the types of the particles and $\left|\bm{\alpha}\right|=\alpha_{1}+\ldots+\alpha_{r}$.

Moreover, we have introduced two classes of models for $r=2$. When each vertex may produce no children or one child of both types, the model is called full binary tree. In addition, if we allow to each vertex to produce one child of the same own type then the model is called full binary tree with survivals. In the study of the likelihood of the multitype Galton-Watson trees concerning $\bm{\mu}$, we show a relation with the \textit{Narayana numbers} $N\left(k,l\right)$ (see for more details \cite{K4}),
$$
N\left(k,l\right)=\frac{1}{k}\binom{k}{l}\binom{k}{l-1},k\geq1,l=1,\ldots,k.
$$
The Narayana numbers are a sort of generalization of the \textit{Catalan numbers} (see \cite{K5}), and from them we have obtained, in a non linear way, the likelihood of the number of the type 1 ($n$) and type 2 ($m$) vertices having exactly two children,  
\begin{equation}\label{eq:distribuione_padri_dx_e_sx}
\mathscr{L}\left(P,Q|n,m\right)=N\left(n+m,m+1\right)\cdot P^{m+1}\left(1-P\right)^{n}Q^{n}\left(1-Q\right)^{m},
\end{equation}
where the paramters $P,Q\in\left(0,1\right)$ depend only on the offspring distribution $\bm{\mu}=\left(\mu^{\left(1\right)},\mu^{\left(2\right)}\right)$ of the tree. Moreover, through the Narayana numbers, we have outlined an interesting rappresentation of the \textit{full binary trees} as two-dimensional decompositions of the integers (see for more details \cite{K6}).

In Section 1 we give the basical definitions of the two-dimensional decompositions of integers and the Narayana numbers and, in a particular case, we prove a connection between them. The full binary trees are briefly defined in Section 2, and we outline a relation between them and the Narayana numbers. In Section 1 we prove the characterization of full binary trees as particular two-dimensional decompositions of integers. Then the preliminaries on multitype Galton-Watson trees are given in Section 3. We illustrate in this Section a characterization of the moment-generating function of the period of the contour function of the trees. Eventually, in Section 4 we give the rigorous definitions of full binary trees and the full binary trees with survivals and we prove a sufficient condition for the period of the contour process to be finite. We show the likelihood \eqref{eq:distribuione_padri_dx_e_sx} as final result.  

\section{Two-dimensional decompositions of integers and Narayana numbers}\label{sec:relazione_decomposizizoni_narayana}
\noindent At first we give the basical definitions of the Narayana numbers (see \cite[Abstract and Section 1.1]{K4}) and the two-dimensional decompositions of integers (see \cite[Vol.2, Section IX, Chapter II, Paragraph 429]{K6}, then, in a particular case, we prove a relation between them.
\begin{defin}[Narayana numbers]\label{def:narayana_numbers}
The Narayana numbers $N\left(n,k\right)$, $n\geq1$ and $k=1,\ldots,n$ are defined in the followinf way
\begin{equation*}\label{eq:naraya_numbers}
N\left(n,k\right)=\frac{1}{n}\binom{n}{k}\binom{n}{k-1},
\end{equation*}
and they are usefull for counting problems. For example from \cite[Section 1.1]{K4}, it is known that the Narayana number $N\left(n,k\right)$ is the number of expressions containing $n$ pairs of parentheses which are correctly matched and which contain $k$ distinct nestings. For instance, $N\left(4,3\right)=6$ counts all the following expressions with $4$ pairs of parentheses, which each contains three times the sub-pattern ( ),
$$
(\ ()()()\ ),\ \ (\ ()()\ )(),\ \ (\ ()\ )()(),\ \ ()(\ ()()\ ),\ \ ()()(\ ()\ ),\ \ ()(\ ()\ )().
$$ 
\end{defin} 
\begin{defin}[Two-dimensional decompostion of integers]\label{def:decomposizione_numeri_interi}
Let $d\geq 1$, $b\geq 1$, $c\geq 0$ and $w\geq 0$ be integers. Consider a matrix $b\times d$ with elements limited in magnitude to $c$ (zero being included) and in descending order in each row and column, and such that the sum of all the elements is exactly $w$ (see Figure \ref{fig:esempio_matrice_di_decomposizione}). Then, each of these matrices is called a \textit{two-dimensional decomposition} of $w$ with parameters $d,b,c$.
\begin{figure}[t]
\centering
$$
\begin{bmatrix}
a_{1,1} & \geq & \cdots & \geq & a_{1,d} \\
\ugec & & \ugec & & \ugec \\
\cdots & & \cdots & & \cdots \\
\ugec & & \ugec & & \ugec \\
a_{b,1} & \geq & \cdots & \geq & a_{b,d}
\end{bmatrix},
$$
with $\forall\ i=1,\ldots,b,\ j=1,\ldots,d,\ a_{i,j}\in\{0,\ldots,c\}$ and $\sum a_{i,j}=w$
\caption{A two-dimensional decomposition with parameters $a,b,c$.}
\label{fig:esempio_matrice_di_decomposizione}
\end{figure}
\end{defin}
Thus, consider the following result about the Narayana numbers and the two-dimensional decompositions of integers.
\begin{lem}\label{lem:relazione_narayana_numeri_e_decomposizioni_interi}
In accordance to Definition \ref{def:decomposizione_numeri_interi}, let $b=2,\ d\geq 1,\ c\geq 0$ and $w\geq 0$ be integers. Then
\begin{equation}\label{eq:relazione_decomposizioni_numeri_narayana}
N\left(c+d+1,d+1\right)=\sum_{w=0}^{2dc}C_{w},
\end{equation}
where $C_{w}$ is the number of the two-dimensional decompositions of $w$.
\end{lem}
\proof
From \cite{K6}, we have that given the following function
\begin{equation}\label{eq:generating_function_two_dim_decomposition_nostro caso}
GF_{d,c}\left(x\right)=\frac{\left(1-x^{c+2}\right)\cdots \left(1-x^{c+d+1}\right)\cdot \left(1-x^{c+1}\right)\cdots \left(1-x^{c+d}\right)}{\left(1-x^{2}\right)\cdots \left(1-x^{d+1}\right)\cdot \left(1-x\right)\cdots \left(1-x^{d}\right)}, 
\end{equation}
the number of two-dimensional decomposition of $w$ is the coefficient $C_{w}$ of $x^{w}$ in $GF_{d,c}\left(x\right)$, written in the power series form. Note that we can only represent integers in the set $\{0,1,\ldots,2dc\}$, and so we have that
\begin{equation}\label{eq:uguaglianza_GF_coefficienti_w}
GF_{d,c}\left(1\right)=\sum_{w=0}^{2dc}C_{w},
\end{equation}
and moreover, using the equality
$$\left(1-x\right)\left(1+x+\cdots+x^{n}\right)=1-x^{n+1},\ \forall n\geq0,$$
the function \eqref{eq:generating_function_two_dim_decomposition_nostro caso} becomes
\begin{equation}\label{eq:semplificazione_GF}
GF_{d,c}\left(x\right)=\displaystyle\frac{\displaystyle\prod_{i=1}^{d}\left(1+\cdots+x^{c+i}\right)\cdot\prod_{i=0}^{d-1}\left(1+\cdots+x^{c+i}\right)}{\displaystyle\prod_{i=1}^{d}\left(1+\cdots+x^{i}\right)\cdot\prod_{i=0}^{d-1}\left(1+\cdots+x^{i}\right)}.
\end{equation}
Now, we compute \eqref{eq:semplificazione_GF} for $x=1$ and get
\begin{align*}
GF_{d,c}\left(1\right) &=\frac{\left(c+2\right)\cdots \left(c+d+1\right)\cdot \left(c+1\right)\cdots \left(c+d\right)}{\left(2\right)\cdots \left(d+1\right)\cdot \left(1\right)\cdots \left(d\right)}\\
& =\frac{\left(c+d+1\right)!\left(c+d\right)!}{\left(c+1\right)!\left(c\right)!\left(d+1\right)!\left(d\right)!}=\binom{c+d+1}{d+1}\frac{\left(c+d\right)!}{\left(c+1\right)!\left(d\right)!}
\\
& =\binom{c+d+1}{d+1}\frac{\left(c+d+1\right)!}{\left(c+1\right)!\left(d\right)!}\cdot\frac{1}{c+d+1}
\\
& = N\left(c+d+1,d+1\right)=\sum_{w=0}^{2dc}C_{w}.
\end{align*}
In this way we can include the particular case of decompositions with $d=0$, indeed it is sufficient to fix $GF_{0,c}\left(1\right)=N\left(c+1,1\right)=1$.
\endproof
\section{Relation between full binary trees and two dimensional decompositions of integers}\label{sec:sezione_2}
\noindent In this section we show an important relation between the full binary trees and two dimensional decompositions of integers. 
\subsection*{Full binary trees} A full binary tree, the rigorous definition of which is given in Section \ref{sec:fulbintee_withsur}, represents a two-type system particles in which every type $i$ particle may produce no particles or exactly two particles (i.e.\ fathers, see Section \ref{sec:distributione_padri_alberi_bin_con_sopr}), one of type 1 and the other one of type 2, for every $i=1,2$. In this Section we consider the special case of finite full binary trees, the condition of a.s.\ finiteness is proved in the Lemma \ref{lem:critical_discussion_binary_model} in Section \ref{sec:distributione_padri_alberi_bin_con_sopr}. In accordance to the Defintion \ref{def:nostro_ordine_tipi_figli} in Section \ref{sec:preliminaries_on_trees} we assume that the left children are considered as vertices of type 1 and the right children are considered as vertices of type 2.
By the definition of the full binary trees and since the root is assume to be a type 1 vertex, it is easy to check that if $n$ and $m$ are the numbers of left and right fathers respectively, then the numbers of the left and right leaves are equal to $m+1$ and $n$ respectively.
\subsection*{Encoding of full binary trees}\label{par:codifica_LR_e_parentesi}  We give the ``()" encoding of the trees, i.e.\ the ``LR" encoding presented by R.P. Grimaldi in \cite[Chapter 24, Example 24.3]{K5} in which we substitute each L with an open parenthesis (, and each R with close one ): if we are in a vertex with children we at first visit the left one, if we are in a left leaf we visit its right brother and if we are in a right leaf we visit the older not visited right vertex having the youngest last ancestor in common, and we write L or R for each left or right vertex visited. Note that each couple of consecutive parentheses ``()" represents a left leaf. Thus, the total number of ( or ), is the number of the fathers and the number of couple ``()" is equal to $m+1$.
\begin{rem}\label{oss:unicita_LR_parentesi_codifica}
Note that each full binary tree has an unique ``()" encoding.
\end{rem}
\begin{figure}[h!]
\centering
\begin{tikzpicture}[scale=0.7,every node/.style={scale=0.7,draw,shape=circle,fill=black},level 1/.style={sibling distance=100mm},level 2/.style={sibling distance=50mm},level 3/.style={sibling distance=30mm},level distance = .9 cm]
    \node (nA)[label=0:$\varnothing$]{}
   child { node (nB) [label=180:{L,(}]{}
              child { node (nD) [label=0:{L,(}]{}
              child { node (nI) [label=270:{L,(}]{} }
                         child { node (nJ) [label=270:{R, )}]{} }
                       }
              child {  node (nE) [label=0:{R, )}]{}
              child { node (nAA) [label=180:{L,(}]{}  }
                         child { node (nBB) [label=180:{R, )}]{} 
                         child { node (nCC) [label=180:{L,(}]{} }
                         child { node (nDD) [label=180:{R, )}]{} }}                         
                       }
            }
   child { node (nC) [label=0:{R, )}]{}
              child { node (nF) [label=0:{L,(}]{}
                         child { node (nK) [label=0:{L,(}]{}  }
                         child { node (nL) [label=0:{R, )}]{}
                         child { node (nEE) [label=0:{L,(}]{}
                         child { node (nGG) [label=270:{L,(}]{}  }
                         child { node (nHH) [label=270:{R, )}]{} }  }
                         child { node (nFF) [label=0:{R, )}]{} } }
                       }
              child {  node (nG) [label=0:{R, )}]{} }
             };
\end{tikzpicture}
\caption{A full-binary tree with $5$ left fathers (including the root) and $4$ right fathers.}  
\label{fig:L_R_encoding}
\end{figure}
\begin{ese}
Let us consider the full binary tree, with a type 1 root, in Figure \ref{fig:L_R_encoding}. It has $5$ left fathers (including the root) and $4$ right fathers. 
The ``()" encoding is (( () ) () () )( () ( () )).
\end{ese}
Thus, we have the following result
\begin{lem}\label{lem:problema_conteggio_padri_con_parentesi}
Let $n\geq 1$ and $m\geq 0$. Then the number of the full binary trees with exactly $n$ left fathers vertices (included the root) and $m$ right fathers vertices is the Narayana number 
$$
N\left(n+m,m+1\right)=\frac{1}{n+m}\binom{n+m}{m+1}\binom{n+m}{m}.
$$
\end{lem}
\proof
It is sufficient to apply the property of the Narayana numbers in Defintion \ref{def:narayana_numbers} in Section \ref{sec:relazione_decomposizizoni_narayana} to the ``()" encoding of the full binary trees seen above.
\endproof  
It is obvious that for every $l\geq 1$ and $k=1,\ldots,l$, the Narayana number $N\left(l,k\right)$ counts the number of full binary trees with  exactly $l-k+1$ left fathers vertices and $k-1$ right fathers vertices. From \eqref{eq:relazione_decomposizioni_numeri_narayana} we have an important relation between the full binary trees and the two-dimensional decompositions, i.e.\ fixed $d\geq 0$ and $c\geq 0$, let be $GF_{d,c}\left(x\right)$ the function in \eqref{eq:generating_function_two_dim_decomposition_nostro caso} and $C_{w}$ the number of the two-dimensional decompositions of the integer $w\in\{0,\ldots,2dc\}$ then, from \eqref{eq:relazione_decomposizioni_numeri_narayana}, we have
\begin{equation}\label{eq:prima_relazione_decomposizioni_narayana}
\sum_{w=0}^{2dc}C_{w}=
\end{equation}
\begin{center}
= number of full binary trees with $\left(c+1\right)$ left fathers (included the root) and $d$ right fathers.
\end{center}
Moreover, extending the result expressed by \eqref{eq:prima_relazione_decomposizioni_narayana}, we are also able to pass from a full binary tree to a two-dimensional decompostion and vice versa.
\begin{teo}
Let $d\geq 1$, $c\geq 0$. Denote with $R_{dc}$ and $FBT_{d,c+1}$ the set of the two-dimensional decompositions with parameters $d,\ c$ of integers $\{0,\ldots,2dc\}$ and the set of the full binary trees with $c+1$ left fathers (including the root) and $d$ right fathers, respectively.\\
Then, exists an unique representation for each element of $R_{dc}$ in $FBT_{d,c+1}$ and vice versa.    
\end{teo}
\proof
To pass from a full binary tree with $c+1$ left fathers and $d$ right fathers to its related two-dimensional decomposition we have to compute the ``LR" and the ``parentheses"" encodings of the tree (see Section \ref{sec:distributione_padri_alberi_bin_con_sopr}), i.e.\ a string of $d+c+1$ couples of parentheses ( , ) with exactly $d+1$ distinct nestings (). Therefore there are $c$ separated couples ( , ) in the string. Then we consider the following defintions: $a_{1,i}$ is the number of ) of the separated couples (,) that stay after the $(i+1)$-$th$ nesting (), for every $i=1,\ldots,d$. If $a_{1,1}<c$ then the remaining $c-a_{1,1}$ ) are all between the first and second nesting (). Similarly, the element $a_{2,i}$ is the number of ( between the separated couples (,) that stay after the $i$-$th$ nesting (), for every $i=1,\ldots,d$ and if $a_{2,1}<c$ then the remaining $c-a_{2,1}$ ( are all before the first nesting (). So, according to the defintions of $a_{i,j}$, the two-dimensional decomposition related to the tree is 
$$
\begin{bmatrix}
a_{1,1} & \geq & \cdots & \geq & a_{1,d} \\
\ugec & & \ugec & & \ugec \\
a_{2,1} & \geq & \cdots & \geq & a_{2,d}
\end{bmatrix}.
$$
Note that the descending order of the rows is verified from the definitions of elements $a_{i,j}$ and, because the number of the open parentheses ( can not be greater than the closed one, even the descending order of the columns is also verified.\\
To pass from a two dimensional decomposition with parameters $d$ and $c$ to its related full binary tree with $c+1$ left fathers (including the root) and $d$ right fathers it is enough to consider the definitions of the elements $a_{i,j}$ given above and write the related string of parentheses which encode a full binary tree with $c+1$ left fathers and $d$ right fathers. The uniqueness of the representations is given respectively by the definitions of the elements of the decompositions and the uniqueness of the ``LR" encoding of the trees (see Remark \ref{oss:unicita_LR_parentesi_codifica}). 
\endproof
\begin{ese}\label{es:esempio_alberi_matrici}
Let $d=2,c=1$, so by \eqref{eq:prima_relazione_decomposizioni_narayana} we have that
$$
\sum_{w=0}^{4}C_{w}=N\left(4,3\right)=6.
$$
\begin{table}[t]
\centering
\begin{footnotesize}
\begin{tabular}{|c|c|c|}
\hline $\begin{bmatrix}
0 & 0 \\ 0 & 0
\end{bmatrix}$,\ $w=0$ & $\begin{bmatrix}
1 & 0 \\ 0 & 0
\end{bmatrix}$,\ $w=1$ & $\begin{bmatrix}
1 & 1 \\ 0 & 0
\end{bmatrix}$,\ $w=2$ \\ 
\hline $\begin{tikzpicture}[scale=0.65, every node/.style={scale=0.65,draw,shape=circle,fill=black},level 1/.style={sibling distance=16mm},level 2/.style={sibling distance=8mm},level 3/.style={sibling distance=8mm},level distance = .8 cm]
    \node (nA)[label=90:$\varnothing$]{}
   child {node (nB)[label=90:L]{}
    child {node (nC)[label=270:L]{}}
    child {node (nD)[label=270:R]{}}}
   child {node (nE)[label=0:R]{}
    child {node (nF)[label=270:L]{}}
    child {node (nG)[label=0:R]{}
     child {node (nH)[label=270:L]{}}
     child {node (nI)[label=270:R]{}}}} ;
\end{tikzpicture}$ & $\begin{tikzpicture}[scale=0.65,every node/.style={scale=0.65,draw,shape=circle,fill=black},level 1/.style={sibling distance=16mm},level 2/.style={sibling distance=8mm},level 3/.style={sibling distance=8mm},level distance = .8 cm]
    \node (nA)[label=90:$\varnothing$]{}
   child {node (nB)[label=180:L]{}
   child {node (nC)[label=270:L]{}}
   child {node (nD)[label=270:R]{}
   child {node (nE)[label=270:L]{}}
   child {node (nF)[label=270:R]{}}}}
   child {node (nG)[label=0:R]{}
   child {node (nH)[label=270:L]{}}
   child {node (nH)[label=270:R]{}}};
\end{tikzpicture}$ & $\begin{tikzpicture}[scale=0.65,every node/.style={scale=0.65,draw,shape=circle,fill=black},level 1/.style={sibling distance=16mm},level 2/.style={sibling distance=8mm},level 3/.style={sibling distance=8mm},level distance = .8 cm]
    \node (nA)[label=0:$\varnothing$]{}
   child {node (nB)[label=180:L]{}
    child {node (nC)[label=127:L]{}}
    child {node (nD)[label=0:R]{}
     child {node (nE)[label=270:L]{}}
     child {node (nF)[label=0:R]{}
      child {node (nG)[label=270:L]{}}
      child {node (nH)[label=270:R]{}}}}}
   child {node (nI)[label=270:R]{}}       
    ;
    \end{tikzpicture}$ \\
\hline $\left(\left(\right)\right)\left(\right)\left(\right)$ & $\left(\left(\right)\left(\right)\right)\left(\right)$ & $\left(\left(\right)\left(\right)\left(\right)\right)$ \\
\hline $\begin{bmatrix}
1 & 1 \\ 1 & 1
\end{bmatrix}$,\ $w=4$ & $\begin{bmatrix}
1 & 0 \\ 1 & 0
\end{bmatrix}$,\ $w=2$ & $\begin{bmatrix}
1 & 1 \\ 1 & 0
\end{bmatrix}$,\ $w=3$ \\
\hline $\begin{tikzpicture}[scale=0.65,every node/.style={scale=0.65,draw,shape=circle,fill=black},level 1/.style={sibling distance=16mm},level 2/.style={sibling distance=8mm},level 3/.style={sibling distance=8mm},level distance = .8 cm]
    \node (nA)[label=0:$\varnothing$]{}
   child {node (nB)[label=0:L]{}}
   child {node (nC)[label=0:R]{}
    child {node (nD)[label=180:L]{}}
    child {node (nE)[label=0:R]{}
     child {node (nF)[label=180:L]{}
      child {node (nG)[label=180:L]{}}
      child {node (nH)[label=0:R]{}}}
     child {node (nI)[label=0:R]{}}}}
   ;
    \end{tikzpicture}$ & $\begin{tikzpicture}[scale=0.65,every node/.style={scale=0.65,draw,shape=circle,fill=black},level 1/.style={sibling distance=16mm},level 2/.style={sibling distance=16mm},level 3/.style={sibling distance=8mm},level distance = .8 cm]
    \node (nA)[label=90:$\varnothing$]{}
  child {node (nB)[label=90:L]{}}
  child {node (nC)[label=0:R]{}
   child {node (nD)[label=180:L]{}
    child {node (nF)[label=270:L]{}}
    child {node (nG)[label=270:R]{}}}
   child {node (nE)[label=180:R]{}
    child {node (nH)[label=270:L]{}}
    child {node (nI)[label=270:R]{}}}} 
   ;
    \end{tikzpicture}$ & $\begin{tikzpicture}[scale=0.65,every node/.style={scale=0.65,draw,shape=circle,fill=black},level 1/.style={sibling distance=16mm},level 2/.style={sibling distance=16mm},level 3/.style={sibling distance=8mm},level distance = .8 cm]
    \node (nA)[label=0:$\varnothing$]{}
  child {node (nB)[label=0:L]{}}
  child {node (nC)[label=0:R]{}
   child {node (nD)[label=180:L]{}
    child {node (nF)[label=180:L]{}}
    child {node (nG)[label=0:R]{}
     child {node (nH)[label=270:L]{}}
     child {node (nI)[label=270:R]{}}}}
   child {node (nE)[label=270:R]{}}} 
   ;
    \end{tikzpicture}$ \\
\hline $\left(\right)\left(\right)\left(\left(\right)\right)$ & $\left(\right)\left(\left(\right)\right)\left(\right)$ & $\left(\right)\left(\left(\right)\left(\right)\right)$ \\
\hline
\end{tabular}
\end{footnotesize}

\hspace{1mm}
\caption{Representations of finite full binary trees with $2$ left fathers (including the root) and $2$ right fathers as two-dimensional decompositions of integers.}
\label{tab:seconda_parte_esempio_rappresentazione_alberi_come_decomposizione}
\end{table}
In Table \ref{tab:seconda_parte_esempio_rappresentazione_alberi_come_decomposizione} you can see the rappresentations of all the 6 full binary trees with 2 left fathers and 2 right fathers, in accordance with the algorithm seen above. For each tree is shown the ``LR" encoding, the related ``parentheses" encoding and the two-dimensional decompostion related and moreover the integer $w$ rappresented in such decomposition.
\end{ese}
\section{Multitype Galton-Watson trees}\label{sec:preliminaries_on_trees}
\noindent We refer to G. Miermont \cite[Sec. 1.3, 1.4]{K1} for the notions of the multitype trees and the multitype Galton-Watson trees. 
\subsection{Multitype trees}
\noindent For $n\geq 0$, let $U$ be the infinite-regular tree
\begin{equation*}
U=\bigcup_{n\geq 0}\mathbb{N}^{n},
\end{equation*}
where if $\textbf{\textit{u}}\in U$, then $\textbf{\textit{u}}=\left(u_{1},\ldots,u_{n}\right)\in\mathbb{N}^{n},\ u_{i}\in\mathbb{N},\ i=1,\ldots,n$. We use the convention $\mathbb{N}^{0}=\{\varnothing\}$ throughout. For $\textbf{\textit{u}}=\left(u_{1},\ldots,u_{n}\right),\textbf{\textit{v}}=\left(v_{1},\ldots,v_{m}\right)\in U$, we let 
$$\textbf{\textit{u}}\textbf{\textit{v}}=\left(u_{1},\ldots,u_{n},v_{1},\ldots,v_{m}\right)\in U$$ be their concatenation and $\left|\textbf{\textit{u}}\right|=n,\left|\textbf{\textit{v}}\right|=m$ their length (with the convention $\left|\varnothing\right|=0$). Let $\textbf{\textit{u}}\in U$ and $A\subseteq U$, we let $\textbf{\textit{u}}A=\{\textbf{\textit{u}}\textbf{\textit{v}}|\textbf{\textit{v}}\in A\}$, and say that \textbf{\textit{u}} is a \textit{prefix} of \textbf{\textit{w}} if $\textbf{\textit{w}}\in \textbf{\textit{u}}U$, and we write $\textbf{\textit{u}}\vdash \textbf{\textit{w}}$.\\
Now we give the definition of a \textit{planar tree}
\begin{defin}[Planar tree]
A planar tree is a finite subset $\tau$ of $U$ such that
\begin{itemize}
\item[$\left(i\right)$] $\varnothing\in \tau$, and it is called the \textit{root} of $\tau$,
\item[$\left(ii\right)$] $\forall\ \textbf{\textit{u}}\in U$ and $i\in \mathbb{N}$, if $\textbf{\textit{u}}i\in \tau\Rightarrow\ \textbf{\textit{u}}\in \tau$, and $\textbf{\textit{u}}j\in \tau$ for every $1\leq j\leq i$. 
\end{itemize}
\end{defin} 
\noindent Moreover, an element $\textbf{\textit{u}}\in \tau$ is called a \textit{vertex} of $\tau$, and $\Vert\tau\Vert$ is the number of edges of the tree $\tau$. We let $T$ be the set of all planar trees, which we refer to as \textit{trees} in the sequel. Now we give some important definitions about trees.
\begin{defin}
Let $\tau\in T$ and \textbf{\textit{u}}, \textbf{\textit{v}} $\in T$ then
\begin{itemize}
\item the number $c_{\tau}\left(\textbf{\textit{u}}\right)=\max\{i\in\mathbb{N}^{+}|\textbf{\textit{u}}i\in \tau \},\text{ with }\textbf{\textit{u}}0=\textbf{\textit{u}}$ is the number of children of \textbf{\textit{u}},
\item the set of the \textit{leaves} of $\tau$ is defined as $\{\textbf{\textit{u}}\in \tau|c_{\tau}\left(\textbf{\textit{u}}\right)=0\}$,
\item \textbf{\textit{u}} is an \textit{ancestor} of \textbf{\textit{v}} if $\textbf{\textit{u}}\vdash\textbf{\textit{v}}$. 
\end{itemize}
\end{defin}
\noindent Any tree $\tau\in T$ is endowed with the \textit{depth-first order},   
\begin{defin}[Depth-first order $\prec$]
Let $\tau\in T$ and $\textbf{\textit{u}},\textbf{\textit{v}}\in \tau$, then
$$
\textbf{\textit{u}}\prec \textbf{\textit{v}}\text{ if }
\textbf{\textit{u}}\vdash \textbf{\textit{v}}\text{ or }
\textbf{\textit{u}}=\textbf{\textit{w}}\textbf{\textit{u}}^{'},\ \textbf{\textit{v}}=\textbf{\textit{w}}\textbf{\textit{v}}^{'}, \text{where } \textit{u}_{1}^{'}<\textit{v}_{1}^{'} 
.
$$ 
\end{defin} 
\begin{ese}\label{ese:esempio_grafico_albero_planar}
The tree $\tau\in T$ in Figure \ref{fig:esempio_funz_contorno} can be written according to the depth-first order in the following way $\tau=\{\varnothing,1,11,12,2,21,211,212,213,$ $3\}$.
\begin{figure}[h]
\centering
{\begin{tikzpicture}[scale=0.59, every node/.style={scale=0.59,draw,shape=circle,fill=black},level 1/.style={sibling distance=35mm},level 2/.style={sibling distance=20mm},level 3/.style={sibling distance=19mm},level distance = .9 cm]
    \node (nA) [label=0:$\varnothing$]{}
      child { node (nB) [label=90:$1$]{}
       child {node (nC) [label=0:$11$]{}}
       child {node (nD) [label=0:$12$]{}}
      }
      child { node (nE) [label=0:$2$]{}
        child {node (nF) [label=0:$21$]{}
         child {node (nG) [label=0:$211$]{}}
         child {node (nH) [label=0:$212$]{}}
         child {node (nI) [label=0:$213$]{}}}     
       }
      child {node (nL) [label=0:$3$]{}}        
       ;
\end{tikzpicture}}\quad\includegraphics[scale=.35]{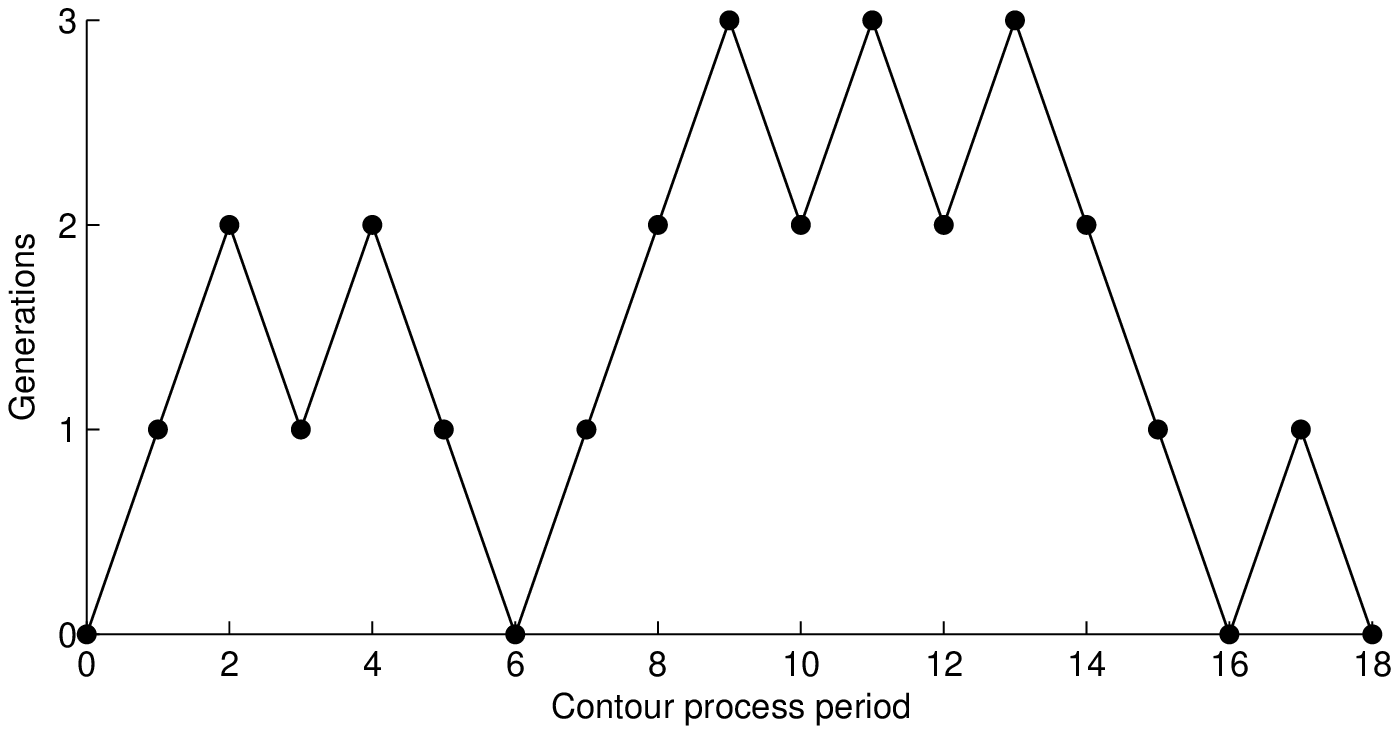}
\caption{A planar tree and its contour process}
\label{fig:esempio_funz_contorno}
\end{figure}\end{ese}
\noindent Now we are able to introduce the $r$-type planar trees, or simply the $r$-type trees.
\begin{defin}[$r$-type planar trees]
Let $r\geq 1$, then a $r$-type planar tree is a pair $\left(\tau,e_{\tau}\right)$, where
\begin{enumerate}
\item $\tau\in T$,
\item $e_{\tau}:\tau\longrightarrow \{1,\ldots,r\}$, i.e. $\forall\ \textbf{\textit{u}}\in \tau,\ e_{\tau}\left(\textbf{\textit{u}}\right)\in\{1,\ldots,r\}$ is called the type of the vertex $\textbf{\textit{u}}$.
\end{enumerate}
Moreover, let $T^{\left(r\right)}$ be the set of $r$-type trees and we let 
$$
T^{\left(r\right)}_{i}=\{\tau\in T^{\left(r\right)}|e_{\tau}\left(\varnothing\right)=i\}\ \forall\ i\in\{1,\ldots,r\}.
$$ 
\end{defin}
\noindent Now the purpose is to count the children of a vertex, according to the type. To do this, at first we define the \textit{counter map}.
\begin{defin}[Counter map]\label{def:counter_map}
Let $r\geq 1$ and $\displaystyle W_{r}=\bigcup_{n\geq 0}\{1,\ldots,r\}^{n}$ be the set of finite, possibly empty, $\{1,\ldots,r\}$-valued sequences, then the counter map, $\textbf{\textit{p}}:W_{r}\longrightarrow \mathbb{N}^{r}$, is such that
$$
\textbf{\textit{p}}\left(\textbf{\textit{w}}\right)=\left(p_{1}\left(\textbf{\textit{w}}\right),\ldots,p_{r}\left(\textbf{\textit{w}}\right)\right)\ \forall\ \textbf{\textit{w}}\in W_{r},
$$
where $p_{i}\left(\textbf{\textit{w}}\right)$ is the number of $i$ in $\textbf{\textit{w}}$, $\forall\ i=1,\ldots,r$. 
\end{defin}
\noindent So, $\forall\ \left(\tau,e_{\tau}\right)\in T^{\left(r\right)}$ and $\forall\ \textbf{\textit{u}}\in \tau$ we can define the following vector,
\begin{equation*}
\textbf{\textit{w}}_{\tau}\left(\textbf{\textit{u}}\right)=\left(e_{\tau}\left(\textbf{\textit{u}}j\right),1\leq j\leq c_{\tau}\left(\textbf{\textit{u}}\right)\right)\in W_{r},
\end{equation*}
and then
\begin{equation*}
\textbf{\textit{p}}\left(\textbf{\textit{w}}_{\tau}\left(\textbf{\textit{u}}\right)\right)\in \mathbb{N}^{r} 
\end{equation*}
is the vector of the number of children of $\textbf{\textit{u}}$ for each type.
\begin{rem}\label{oss:funz_contorno_alberi_multi}
Note that the graphical representation of a multitype tree is the same of a one-type planar tree. Indeed, neglecting the type of the particles, it may be considered as one-type tree. Moreover, the set of the vertices of a multitype tree is ordered according to the first-depth order.
\end{rem}
\subsection{Galton-Watson trees}\label{sec:galton_watson_trees}
\noindent In this section we treat the multitype planar trees where each vertex has a number of children of certain type according to the offspring distribution.
\begin{defin}[Offspring distribution]\label{def:offspring:distribution}
Let $r\geq 1$ and $\bm{\xi}=\left(\xi^{\left(1\right)},\ldots,\xi^{\left(r\right)}\right)$ be a family of probabilities on the $\sigma$-algebra $\sigma\left(W_{r}\right)$. We say that the family of probabilities $\bm{\mu}$ on the $\sigma$-algebra $\sigma\left(\mathbb{N}^{r}\right)$ is an offspring distribution, where $\textbf{\textit{p}}$ is the counter map defined in Definition \ref{def:counter_map} and 
$$
\bm{\mu}^{\left(i\right)}=\xi^{\left(i\right)}\circ \textbf{\textit{p}}^{-1},\ \ \forall i=1,\ldots,r. 
$$ 
\end{defin}
\noindent Now, we build a distribution on $T^{\left(r\right)}_{i},\ \forall\ i\in\{1,\ldots,r\}$, such that 
\begin{enumerate}
\item different vertices have indipendent offspring
\item type $j$ vertices have a set of children with types given by a sequence $\textbf{\textit{w}}\in W_{r}$ with probability $\xi^{\left(j\right)}\left(\textbf{\textit{w}}\right)$.
\end{enumerate}
To do this, $\forall\ \textbf{\textit{u}}\in U$, let $\textbf{\textit{C}}_{\textbf{\textit{u}}}=\left(C_{\textbf{\textit{u}}}\left(l\right),1\leq l\leq c_{\tau}\left(\textbf{\textit{u}}\right) \right)$ be the vector of the types of the children of $\textbf{\textit{u}}$. $\left(\textbf{\textit{C}}_{\textbf{\textit{u}}}\right)_{\textbf{\textit{u}}\in U}$ is a family of independet vectors and such that $\textbf{\textit{C}}_{\textbf{\textit{u}}}$ has law $\xi^{\left(e_{\tau}\left(\textbf{\textit{u}}\right)\right)}$. Now, recursively, we construct a subset $\tau\subset U$ and a mark-map $e_{\tau}:\tau\rightarrow \{1,\ldots,r\}$ in the following way
\begin{enumerate}
\item $\varnothing \in \tau$
\item $e_{\tau}\left(\varnothing\right)=i$
\item if $\textbf{\textit{u}}\in \tau$, $e\left(\textbf{\textit{u}}\right)=j$, then, with probability $\xi^{\left(j\right)}\left(\textbf{\textit{C}}_{\textbf{\textit{u}}}\right)$, $\textbf{\textit{u}}l\in \tau$ if and only if $1\leq l\leq c_{\tau}\left(\textbf{\textit{u}}\right)$ and then $e\left(\textbf{\textit{u}}l\right)=C_{\textbf{\textit{u}}}\left(l\right)$.
\end{enumerate}  
A pair $\left(\tau,e_{\tau}\right)\in T^{\left(r\right)}_{i}$ equipped with $\bm{\mu}$, for every $r\geq 1$ and $i\in\{1,\ldots,r\}$, is called $\bm{\mu}$-\textit{GW tree}.
\begin{rem}\label{oss:processoGW_da_albero}
It is easy to check that the subset $\tau\subset U$ has the properties of a planar tree (it might be infinite). Moreover, from the construction we have that
\begin{equation}\label{eq:processo_GW_relativo_albero_GW}
\textbf{\textit{Z}}_{n}\left(\tau\right)=\left(\#\{\textbf{\textit{u}}\in \tau:\left|\textbf{\textit{u}}\right|=n,e\left(\textbf{\textit{u}}\right)=i\},\ i\in\{1,\ldots,r\}\right),\ n\geq 0,
\end{equation}
is a multitype Galton-Watson process (with discrete time) with offspring distribution $\bm{\mu}^{\left(i\right)}=\xi^{\left(i\right)}\circ \textbf{\textit{p}}^{-1}$.
\end{rem}
\subsection{Characterization of the period of the contour process of the $\bm{\mu}$-GW trees}
\noindent Let $\left(\tau,e_{\tau}\right)\in T^{\left(r\right)}_{i}$ a planar tree and consider the following order for the offspring of a vertex, according to the type.
\begin{defin}[Type-ordering offspring]\label{def:nostro_ordine_tipi_figli}
Let $\left(\tau,e_{\tau}\right)\in T^{\left(r\right)}_{i}$. Then
\begin{equation}\label{eq:ordinamento_offspring_in_base_al_tipo}
\begin{cases}
e\left(\textbf{\textit{u}}l\right)=1,\ 1\leq l\leq p_{1}\left(\textbf{\textit{w}}_{\tau}\left(\textbf{\textit{u}}\right)\right)\\
\cdots \\
\displaystyle e\left(\textbf{\textit{u}}l\right)=k,\ \sum_{i=1}^{k-1}p_{i}\left(\textbf{\textit{w}}_{\tau}\left(\textbf{\textit{u}}\right)\right)+1\leq l\leq \sum_{i=1}^{k}p_{i}\left(\textbf{\textit{w}}_{\tau}\left(\textbf{\textit{u}}\right)\right)\\
\cdots \\
\displaystyle e\left(\textbf{\textit{u}}l\right)=r,\ \sum_{i=1}^{r-1}p_{i}\left(\textbf{\textit{w}}_{\tau}\left(\textbf{\textit{u}}\right)\right)+1\leq l\leq c_{\tau}\left(\textbf{\textit{u}}\right),
\end{cases}
\end{equation}
where $\textbf{\textit{p}}=\left(p_{1},\ldots,p_{r}\right)$ is the counter map given in Definition \ref{def:counter_map} and $\textbf{\textit{u}}\in \tau$. Note that if $p_{j}\left(\textbf{\textit{w}}_{\tau}\left(\textbf{\textit{u}}\right)\right)=0$ for some $j\in\{1,\ldots,r\}$, then $\{\textbf{\textit{u}}l\in\tau|e\left(\textbf{\textit{u}}l\right)=j\}=\{\varnothing\}$.
\end{defin}    
\noindent From now on, a $\bm{\mu}$-GW tree is equipped with the type-ordering offspring. We can derive the \textit{contour process} 
$$\mathscr{C}_{\tau}:\left[0,2\Vert\tau\Vert\right]\longrightarrow\mathbb{N}$$ of $\left(\tau,e_{\tau}\right)\in T^{\left(r\right)}_{i}$ (see Figure \ref{fig:esempio_funz_contorno}) and \cite[Sec.\ 1.1]{K2} for details. The value $\mathscr{C}_{\tau}\left(n\right)$ at the time $n$ is the generation of the vertex visited at the step $n$ in this evolution.
\noindent Note that $\Vert\tau\Vert$ is a non negative-integer-valued random variable.
\begin{defin}\label{def:numero_edges_nuova_notazione}
Let $\left(\tau,e_{\tau}\right)\in T^{\left(r\right)}$ be a $\bm{\mu}$-GW tree. We define
$
X_{\textbf{\textit{u}},j}
$
the number of edges in the tree $\tau$ from the vertex $\textbf{\textit{u}}$ of type $j$. With this definition we have
$$
X_{\varnothing,j}=\Vert\tau\Vert,
$$
where $\tau$ is a $\bm{\mu}$-GW tree rooted in a type $j$ vertex.
\end{defin}
Recall that $\textbf{\textit{C}}_{\textbf{\textit{u}}}$ is the vector of the types of the children of $\textbf{\textit{u}}$ and $\textbf{\textit{p}}$ is the counter map seen in Defintion \ref{def:counter_map}, and denote $
\textbf{\textit{Y}}_{\textbf{\textit{u}}}=\textbf{\textit{p}}\left(\textbf{\textit{C}}_{\textbf{\textit{u}}}\right)$. Each component
$$
Y_{\textbf{\textit{u}},k}=p_{k}\left(\textbf{\textit{C}}_{\textbf{\textit{u}}}\right),\ k=1,\ldots,r,
$$
is the number of chidren of $\textbf{\textit{u}}$ of type $k$.
The following theorem show some properties of $X_{\textbf{\textit{u}},j}$.
\begin{teo}[Number of edges from a vertex of a $\bm{\mu}$-GW tree]\label{teo:number_of_edges}
Let $\left(\tau,e_{\tau}\right)\in T^{\left(r\right)}$ be a $\bm{\mu}$-GW tree. Then, for every $\textbf{\textit{u}}\in \tau$, we have 
\begin{enumerate}
\item for every $i,j\in\{1,\ldots,r\}$ with $i\neq j$, $X_{\textbf{\textit{u}}l,i},X_{\textbf{\textit{u}}m,j}$ are independent for every $l=1,\ldots,Y_{\textbf{\textit{u}},i}$ and $m=1,\ldots,Y_{\textbf{\textit{u}},i}$,
\item for every $i\in\{1,\ldots,r\}$, $X_{\textbf{\textit{u}}l,i}$ are i.i.d.\ and $X_{\textbf{\textit{u}}l,i}\stackrel{d}{=}X_{\varnothing,i}$ for every $l=1,\ldots,Y_{\textbf{\textit{u}},i}$,
\item 
for every $j\in\{1,\ldots,r\}$ and $\textbf{\textit{u}}\in \tau$, the following recursive formula holds
\begin{equation}\label{eq:definizione_ricorsiva_numero_edges}
X_{\textbf{\textit{u}},j}=\sum_{k=1}^{r}\Big(Y_{\textbf{\textit{u}},k}+\sum_{l=1}^{Y_{\textbf{\textit{u}},k}}X_{\textbf{\textit{u}}l,k}\Big).
\end{equation} 
\end{enumerate}
\end{teo}
\proof
Note that, in accordance with the literature of the Galton-Watson processes, 1 and 2 are verified. In other words, it represents the property that different particles have independent offspring in the Galton-Watson prodesses. 3 can be proved for recursion on the vertices and on the vertices type.
\endproof
As consequence of Theorem 3.1 we can get the characterization of the moment-generating function of the period of $\mathscr{C}_{\tau}$, $2X_{\varnothing,i}$.
\begin{teo}\label{teo:caratterizzazione_funz_gen_momenti_durata}
Let $\left(\tau,e_{\tau}\right)\in T^{\left(r\right)}$, and $F_{i}\left(s\right)=\mathbb{E}\left[e^{s\cdot 2X_{\varnothing,i}}\right]$ be the moment-generating function of $2X_{\varnothing,i}$, $s\in\mathbb{R}$. Then
\begin{equation*}\label{eq:caratt_funz_gen_mom_durata}
F_{i}(s)=\sum_{\bm{\alpha}\in\mathbb{N}^{r}}\Big(\mu^{\left(i\right)}\left(\bm{\alpha}\right)\cdot e^{2s\left|\bm{\alpha}\right|}\cdot\prod_{k=1}^{r}F_{k}(s)^{\alpha_{k}}\Big),\text{ with }\big|\bm{\alpha}\big|=\sum_{i=1}^{r}\alpha_{i}.
\end{equation*}
\end{teo}
\proof Let $f_{k}(s)=\mathbb{E}[e^{s\cdot X_{\varnothing,i}}]$. For the law of total probability we have,
\begin{align*}
f_{i}(s) = &\sum_{\bm{\alpha}\in\mathbb{N}^{r}}\mathbb{E}\Big[\exp(s\cdot X_{\varnothing,i})|\textbf{\textit{Y}}_{\varnothing}=\bm{\alpha}\Big]\cdot \mu^{(i)}(\bm{\alpha})
\\
\stackrel{\eqref{eq:definizione_ricorsiva_numero_edges}}{=}
& 
\sum_{\bm{\alpha}\in\mathbb{N}^{r}}\mathbb{E}\Big[\exp\Big(s\cdot \sum_{k=1}^{r}(\alpha_{k}+\sum_{l=1}^{\alpha_{k}}X_{l,k})\Big)\Big]\cdot \mu^{(i)}(\bm{\alpha})
\\ = &
\sum_{\bm{\alpha}\in\mathbb{N}^{r}}\mu^{(i)}(\bm{\alpha})\cdot \mathbb{E}\Big[\prod_{k=1}^{r}\exp\Big(s\cdot (\alpha_{k}+\sum_{l=1}^{\alpha_{k}}X_{l,k})\Big)\Big]
\intertext{and using 1 and 2 of Theorem \ref{teo:number_of_edges} we get }
f_{i}(s) = &
\sum_{\bm{\alpha}\in\mathbb{N}^{r}}\mu^{(i)}(\bm{\alpha})\cdot \prod_{k=1}^{r}\Big(\exp(s\cdot\alpha_{k})\mathbb{E}\Big[\prod_{l=1}^{\alpha_{k}}\exp(s\cdot X_{l,k})\Big]\Big)
\\
= &
\sum_{\bm{\alpha}\in\mathbb{N}^{r}}\Big(\mu^{(i)}(\bm{\alpha})\cdot e^{s|\bm{\alpha}|}\cdot\prod_{k=1}^{r}f_{k}(s)^{\alpha_{k}}\Big),and \intertext{and hence}
f_{i}(2s)=
&
F_{i}(s)=\sum_{\bm{\alpha}\in\mathbb{N}^{r}}\Big(\mu^{(i)}(\bm{\alpha})\cdot e^{2s|\bm{\alpha}|}\cdot\prod_{k=1}^{r}F_{k}(s)^{\alpha_{k}}\Big). \qedhere
\end{align*}
\endproof

As in the classical case, Theorem \ref{teo:caratterizzazione_funz_gen_momenti_durata} provides a property of the moment-generating function $F_{i}\left(s\right),\forall i\in\{1,\ldots,r\}$.
\begin{cor}\label{cor:valore_funz_generatrice_momenti_in_zero_meno}
Define $F_{i}\left(0^{-}\right):=\displaystyle\lim_{s\rightarrow 0^{-}}F_{i}\left(s\right)$, $\forall i\in\{1,\ldots,r\}$, then
\begin{equation*}
F_{i}\left(0^{-}\right)=\mathbb{P}\left(2X_{\varnothing,i}<\infty\right).
\end{equation*}
\end{cor} 
\proof
Let $i\in\{1,\ldots,r\}$ and $\{\left(p^{\left(i\right)}_{2n}\right)_{n\geq0}\cup\ p^{\left(i\right)}_{\infty}\}$ be the distribution of probability of $2X_{\varnothing,i}$, that is $p^{\left(i\right)}_{2n}=\mathbb{P}\left(2X_{\varnothing,i}=2n\right),n\geq0,$ and $p^{\left(i\right)}_{\infty}=\mathbb{P}\left(2X_{\varnothing,i}=\infty\right)$, with $\displaystyle p^{\left(i\right)}_{\infty}+\sum_{n\geq0}p^{\left(i\right)}_{2n}=1$. Then 
$$
F_{i}\left(s\right)=\mathbb{E}\left[e^{s\cdot2X_{\varnothing,i}}\right]=p^{\left(i\right)}_{\infty}\cdot e^{s\cdot\infty}+\sum_{n\geq0}p^{\left(i\right)}_{2n}\cdot e^{2sn},\ s\in\mathbb{R}.
$$
So, $\forall s<0$ we have
\begin{equation}\label{eq:moment_generating_function_in_s_min_0}
F_{i}\left(s\right)=\sum_{n\geq0}p^{\left(i\right)}_{2n}\cdot e^{2sn},
\end{equation}
and $\forall n\geq0$
$$
\left|e^{2sn}\right|=e^{2sn}\leq 1.
$$
Then, from the Dominated Convergence Theorem, we obtain
$$
\lim_{s\rightarrow0^{-}}F_{i}\left(s\right)=\sum_{n\geq0}p^{\left(i\right)}_{2n}=\mathbb{P}\left(2X_{\varnothing,i}<\infty\right).
\qedhere
$$ 
\endproof
Another consequence is a connection between $F_{i}\left(s\right),\forall i\in\{1,\ldots,r\}$, and the extinction probability of the multitype Galton-Watson processes.
\begin{cor}\label{cor:legame_funz_generatrice_mom_prob_estinzione}
Denote $\textbf{\textit{F}}\left(0^{-}\right)=\left(\mathbb{P}\left(2X_{\varnothing,1}<\infty\right),\ldots,\mathbb{P}\left(2X_{\varnothing,r}<\infty\right)\right)$, then
$$
\textbf{\textit{F}}\left(0^{-}\right)=\textbf{\textit{f}}\left(\textbf{\textit{F}}\left(0^{-}\right)\right),$$
where $\textbf{\textit{f}}\left(\textbf{\textit{s}}\right)$ is the multitype generating function related to the offspring distribution $\bm{\mu}$ (see \cite[Chapter V]{K3}). Moreover, $\textbf{\textit{F}}\left(0^{-}\right)=\textbf{\textit{q}}$.
\end{cor}
\proof
Note that $\forall i\in\{1,\ldots,r\},\forall s<0$ and from \eqref{eq:moment_generating_function_in_s_min_0} we have that
$$
\left|F_{i}\left(s\right)\right|=F_{i}\left(s\right)=\sum_{n\geq0}p^{\left(i\right)}_{2n}\cdot e^{2sn}\leq \sum_{n\geq0}p^{\left(i\right)}_{2n}\leq 1,
$$
and so
$$
\Big|\mu^{(i)}(\bm{\alpha})\cdot e^{2s|\bm{\alpha}|}\cdot\prod_{k=1}^{r}F_{k}(s)^{\alpha_{k}}\Big|=\mu^{(i)}(\bm{\alpha})\cdot e^{2s|\bm{\alpha}|}\cdot\prod_{k=1}^{r}F_{k}(s)^{\alpha_{k}}\leq 1. 
$$
Thus, from the Theorem \ref{teo:caratterizzazione_funz_gen_momenti_durata} and the Theorem of Dominated Convergnce, we have that
$$
F_{i}(0^{-})=\sum_{\bm{\alpha}\in\mathbb{N}^{r}}
\Big(\mu^{(i)}(\bm{\alpha})\cdot\prod_{k=1}^{r}F_{k}(0^{-})^{\alpha_{k}}\Big),\qquad \forall i\in\{1,\ldots,r\},
$$
and so $\textbf{\textit{F}}\left(0^{-}\right)=\textbf{\textit{f}}\left(\textbf{\textit{F}}\left(0^{-}\right)\right)$.
To proof the second part of the corollary it's enough to note that $\forall i\in\{1,\ldots,r\}$,$$2X_{\varnothing,i}<\infty\Longleftrightarrow \textbf{\textit{Z}}_{n}\left(\tau\right)=\textbf{\textit{0}},\ \text{for some }n\geq0,$$ where $\left(\tau,e_{\tau}\right)\in T^{\left(r\right)}_{i}$ is a $\bm{\mu}$-GW tree, and $\textbf{\textit{Z}}_{n}\left(\tau\right)$ is the Galton-Watson process related to the tree \eqref{eq:processo_GW_relativo_albero_GW}. Then, by Corollary \ref{cor:valore_funz_generatrice_momenti_in_zero_meno} and the defintion of the extinction probability $\textbf{\textit{q}}$, we get
$$
F_{i}\left(0^{-}\right)=\mathbb{P}\left(2X_{\varnothing,i}<\infty\right)=\mathbb{P}\left(\textbf{\textit{Z}}_{n}\left(\tau\right)=\textbf{\textit{0}},\ \text{for some }n\geq0\right)=q^{\left(i\right)},
$$ 
and so $\textbf{\textit{F}}\left(0^{-}\right)=\textbf{\textit{q}}$.
\endproof
\section{Full binary trees, full binary trees with survivals and the Narayana numbers}\label{sec:fulbintee_withsur}
\noindent We now characterize the offspring distribution of the the full binary trees with survivals.
\begin{defin}
A \emph{full binary tree with survivals} 
is a $\bm{\mu}$-GW tree $\left(\tau,e_{\tau}\right)\in T^{\left(2\right)}_{i},i=1,2$, equipped with the type-ordering offspring (see Definition \ref{def:nostro_ordine_tipi_figli}), and such that the offspring distribution $\bm{\mu}$ is the following
\begin{equation*}
\begin{cases}
\mu^{\left(1\right)}\left(0,0\right)=p_{0}\\
\mu^{\left(1\right)}\left(1,0\right)=p_{1}\\
\mu^{\left(1\right)}\left(1,1\right)=p_{2}\\
\end{cases},
\begin{cases}
\mu^{\left(2\right)}\left(0,0\right)=q_{0}\\
\mu^{\left(2\right)}\left(0,1\right)=q_{1}\\
\mu^{\left(2\right)}\left(1,1\right)=q_{2}\\
\end{cases}
\end{equation*}
where $p_{i},q_{i}\in\left(0,1\right),i=0,1,2$ and $\sum p_{i}=\sum q_{i}=1$.
\end{defin}
The vertices which produce only one vertex is called the \textit{survivals} and the vertices which produced two vertices are called \textit{fathers}.\\
When  
$
\mu^{\left(1\right)}\left(1,0\right)=\mu^{\left(2\right)}\left(0,1\right)=0,
$
no survivals are expressed in the trees and the model is called \emph{full binary trees without survivals}, or, simply \emph{full binary trees}.
\begin{rem}\label{oss:ordinamento_figli_albero_binarrio_con_nostro_ordine}
Note that, according to the type-ordering offspring defined in \eqref{eq:ordinamento_offspring_in_base_al_tipo} in Section 2, when a vertex of type $i$ produces two vertices then the first one (the left one) is a type 1 vertex and the second one (the right one) is a type 2 vertex. This holds both for the full binary trees and for the full binary trees with survivals.  
\end{rem}
\begin{rem}\label{oss:root_destra_o_sinistro}
From now on, if a full binary tree with survivals (or without survivals) starts with a type $j$ vertex, then the root is considered a left vertex or a right vertex respectively if $j=1$ or $j=2$. Morevoer, in the following sections we consider trees starting with a type 1 vertex (the same arguments are verified with a type 2 root).
\end{rem}
\subsection{Number of fathers of full binary trees with survivals}\label{sec:distributione_padri_alberi_bin_con_sopr}
Now, we consider the full binary trees with survivals. Every particles may live for an unit time interval (with probabilities $p_{1}$, $q_{1}$) and then it may \textit{die} (with probabilities $p_{0}$, $q_{0}$) or live for another unit time interval ($p_{1}$, $q_{1}$) after producing a new particle of the other type (with probabilities $p_{2}$, $q_{2}$). So, at the extinction we can only count how many particles of each type have been produced. Our purpose is to find the likelihood of the number of type 1 and type 2 fathers, conditioning to the a.s.\ finiteness of the tree. In the next Corollary we give the condition for the finiteness of full binary trees with survivals.

\begin{lem}\label{lem:critical_discussion_binary_model}
Let $\left(\tau,e_{\tau}\right)\in T^{\left(2\right)}$ be a full binary tree with survivals. If
$
p_{0}q_{0}-p_{2}q_{2}\geq0,$
then  
$\mathbb{P}\left(\|\tau\|<\infty\right)=1 .$
If
$
p_{0}q_{0}-p_{2}q_{2}<0
$,  then  
$\mathbb{P}\left(\|\tau\|<\infty\right)<1.$
\end{lem}
\proof From Theorem \ref{teo:caratterizzazione_funz_gen_momenti_durata} we have that 
$$\begin{cases}
F_{1}\left(s\right)=\mu^{\left(1\right)}\left(0,0\right)+\mu^{\left(1\right)}\left(1,0\right)e^{2s}F_{1}\left(s\right)+\mu^{\left(1\right)}\left(1,1\right)e^{4s}F_{1}\left(s\right)F_{2}\left(s\right)\\
F_{2}\left(s\right)=\mu^{\left(2\right)}\left(0,0\right)+\mu^{\left(2\right)}\left(0,1\right)e^{2s}F_{2}\left(s\right)+\mu^{\left(2\right)}\left(1,1\right)e^{4s}F_{1}\left(s\right)F_{2}\left(s\right),
\end{cases}$$
and so, passing to the limit $s\rightarrow 0^{-}$, we obtain
$$
\begin{cases}
F_{1}\left(0^{-}\right)=p_{0}+p_{1}F_{1}\left(0^{-}\right)+p_{2}F_{1}\left(0^{-}\right)F_{2}\left(0^{-}\right)\\
F_{2}\left(0^{-}\right)=q_{0}+q_{1}F_{2}\left(0^{-}\right)+q_{2}F_{1}\left(0^{-}\right)F_{2}\left(0^{-}\right).
\end{cases}
$$
It is easy to check that
\begin{itemize}
\item if $p_{0}q_{0}-p_{2}q_{2}\geq 0$ then
$$
\begin{cases}
\displaystyle\frac{p_{0}\left(q_{0}+q_{2}\right)}{q_{2}\left(p_{0}+p_{2}\right)}\geq \frac{p_{2}q_{2}+p_{0}q_{2}}{q_{2}\left(p_{0}+p_{2}\right)}=1\\[4mm]
\displaystyle\frac{q_{0}\left(p_{0}+p_{2}\right)}{p_{2}\left(q_{0}+q_{2}\right)}\geq \frac{p_{2}q_{2}+q_{0}p_{2}}{p_{2}\left(q_{0}+q_{2}\right)}=1,
\end{cases}
$$
and, by Corollary \ref{cor:valore_funz_generatrice_momenti_in_zero_meno}, we have that 
$$
\textbf{\textit{F}}\left(0^{-}\right)=\Big(\mathbb{P}(2L^{\left(1\right)}\left(\varnothing\right)<\infty),\mathbb{P}(2L^{\left(2\right)}\left(\varnothing\right)<\infty)\Big)=(1,1).
$$
\item if $p_{0}q_{0}-p_{2}q_{2}<0$ then from the  Corollary \ref{cor:legame_funz_generatrice_mom_prob_estinzione}; we can use \cite[Theorem~2, pag.~186]{K3} obtaining that
$$
\textbf{\textit{F}}\left(0^{-}\right)=\left(\frac{p_{0}\left(q_{0}+q_{2}\right)}{q_{2}\left(p_{0}+p_{2}\right)},\frac{q_{0}\left(p_{0}+p_{2}\right)}{p_{2}\left(q_{0}+q_{2}\right)}\right)<\left(1,1\right). \qedhere
$$ 
\end{itemize}
\endproof
\begin{rem}\label{oss:condizione_finitezza_albero_full_binary}
For the full binary trees without survivals (i.e., $p_1=q_1=0$) we have that 
$$
p_{0}q_{0}- p_{2}q_{2} \geq 0 \iff p_{2}+q_{2}\leq 1
$$
Indeed
$$
p_{0}q_{0}\geq p_{2}q_{2}\iff\left(1-p_{2}\right)\left(1-q_{2}\right)\geq p_{2}q_{2}
\iff p_{2}+q_{2}\leq 1.
$$
\end{rem}
\subsection*{Likelihood of the number of fathers of full binary trees}
Suppose that $\left(\tau,e_{\tau}\right)\in T^{\left(2\right)}_{1}$ is a full binary tree with survivals with $\mathbb{P}\left(\|\tau\|<\infty\right)=1$. The root is assumed to be a type 1 vertex and we use the following notations
\begin{itemize}
\item
$D_{1}$ is number of type 1 fathers in $\tau$;
\item
$D_{2}$ is number of type 2 fathers in $\tau$;
\item
$S_{1}$ is number of type 1 survivals in $\tau$;
\item
$S_{2}$ is number of type 2 survivals in $\tau$.
\end{itemize}
At first, we want to compute the joint distribution of the number of fathers and survivals for each type, 
\begin{equation}\label{eq:prob_numero_padri_con_sopr}
\mathbb{P}\left(D_{1}=n,D_{2}=m,S_{1}=s_{1},S_{2}=s_{2}\right)
\end{equation}
where $n\geq 1$ and $m,\ s_{1},\ s_{2}\geq 0$.\\ 
The case $n=0$ is not included, see Remark \ref{oss:root_destra_o_sinistro}. To compute the probability \eqref{eq:prob_numero_padri_con_sopr},
we note that the type 1 leaves is equal to $m+1$ and the number of the type 2 leaves is equal to $n$. 
From the ``()" encoding and Lemma \ref{lem:problema_conteggio_padri_con_parentesi} in Section \ref{sec:sezione_2}, 
the number of full binary trees having $n$ left fathers and $m$ right fathers is  $N\left(n+m,m+1\right)$.
Now we insert the survivals, i.e.\ $s_{1},s_{2}\geq 0$. We can choose the $s_{1}$ survivals type 1 vertices 
among $m+n+s_{1}$ type 1 vertices of the tree, excluded the last one that is certainly a leaf. 
For the same reason, $s_{2}$ can be chosen among $m+n+s_{2}-1$ vertices. 
Thus, the probability \eqref{eq:prob_numero_padri_con_sopr} becomes
\begin{align*} 
\mathbb{P}(D_{1}=n,D_{2}=m,& S_{1}=s_{1},S_{2}=s_{2}) \\
\notag
& =N\left(n+m,m+1\right)\binom{m+n+s_{1}}{s_{1}}\binom{m+n+s_{2}-1}{s_{2}}\\
\notag
& \qquad\qquad p_{0}^{m+1}\cdot  p_{1}^{s_{1}}\cdot p_{2}^{n}\cdot q_{0}^{n}\cdot q_{1}^{s_{2}}\cdot q_{2}^{m}
\end{align*} 
where $\ n\geq 1$ and $m,s_{1},s_{2}\geq 0$.
The negative binomial distribution implies that \\
$\sum_{s_{1}\geq 0}\binom{m+n+s_{1}}{s_{1}}p_{1}^{s_{1}}\left(1-p_{1}\right)^{n+m+1}=1$
and $\sum_{s_{2}\geq 0}\binom{m+n+s_{2}-1}{s_{2}}q_{1}^{s_{2}}\left(1-q_{1}\right)^{n+m}=1$.
Hence 
\begin{equation}\label{eq:probabilita_congiunta}
\mathbb{P}\left(D_{1}=n,D_{2}=m\right)=N\left(n+m,m+1\right)\cdot\frac{p_{0}^{m+1}p_{2}^{n}}{\left(p_{0}+p_{2}\right)^{n+m+1}}\cdot\frac{q_{0}^{n}q_{2}^{m}}{\left(q_{0}+q_{2}\right)^{n+m}}.
\end{equation}
If we denote
\begin{equation*}
P=\frac{p_{0}}{p_{0}+p_{2}}\text{  and  }Q=\frac{q_{0}}{q_{0}+q_{2}},
\end{equation*}
we get the likelihood of the fathers of the full binary trees (with survivals).
\begin{teo}
Let $\left(\tau, e_{\tau}\right)\in T^{\left(2\right)}_{1}$ be a full binary tree 
with survivals such that $\mathbb{P}\left(\|\tau\|<\infty\right)=1$. 
It $\tau$ has $n$ left fathers and $m$ right fathers, with $n\geq 1$ and $m\geq 0$, then the likelihood reads
\begin{equation*}\label{eq:probl_congiunta_semplificata}
\mathscr{L}\left(P,Q|n,m\right)=N\left(n+m,m+1\right)\cdot P^{m+1}\left(1-P\right)^{n}Q^{n}\left(1-Q\right)^{m},
\end{equation*}
where parameters $P,Q\in\left(0,1\right)$ and depending only from the offspring distribution $\bm{\mu}$ of the tree.
\end{teo}
As a consequence, from the theory of the likelihood estimators (see for instance \cite{K7}), 
the estimators of the parameters $P$, $Q$ are 
\begin{equation*}\label{eq:estimatori_massimi_distribuzione_padri_con_sopr}
\begin{cases}
\displaystyle \widetilde{P}=\frac{m+1}{m+n+1}\\[4mm]
\displaystyle \widetilde{Q}=\frac{n}{m+n}
\end{cases},\text{ where }\tilde{P},\tilde{Q}\in\left(0,1\right),
\end{equation*}
and hence
$$
\Big(\widetilde{\frac{p_{2}}{p_{0}}} \Big)=\frac{n}{m+1}\qquad \text{and}\qquad \Big(\widetilde{\frac{q_{2}}{q_{0}}}\Big)=\frac{m}{n}. 
$$

\end{document}